%\def\Bbb{\bf}

% This file defines the standard font setup, including the macros
% \tenpoint and \twelvepoint, as well as some standard dimension settings.
%
% It is indended to be used as a basis for most other ``personal'' formats.

% Make sure we haven't already been run

\ifx\UsualIsLoaded\undefined
\let\UsualIsLoaded=\relax		% define UsualIsLoaded

\font\fourteenrm=cmr12  scaled \magstep1
\font\fourteenbf=cmbx12 scaled \magstep1
\font\fourteentt=cmtt12 scaled \magstep1
\font\fourteensl=cmsl12 scaled \magstep1
\font\fourteensy=cmsy10 scaled \magstep2
\font\fourteeni=cmmi12  scaled \magstep1
\font\fourteenit=cmti12 scaled \magstep1
\font\fourteensc=cmcsc10 scaled \magstep2
\font\fourteenbit=cmssi12 scaled \magstep1
\font\fourteenbbb=msbm10 scaled \magstep2

\font\twelverm=cmr12
\font\twelvebf=cmbx12
\font\twelvett=cmtt12
\font\twelvesl=cmsl12
\font\twelvesy=cmsy10 scaled \magstep1
\font\twelvei=cmmi12
\font\twelveit=cmti12
\font\twelvesc=cmcsc10 scaled \magstep1
\font\twelvebit=cmssi12
\font\twelvebbb=msbm10 scaled \magstep1
\font\tenrm=cmr10
\font\tenbf=cmb10 
\font\tentt=cmtt10 
\font\tensl=cmsl10 
\font\tensy=cmsy10 
\font\teni=cmmi10 
\font\tenit=cmti10 
\font\tensc=cmcsc10 
\font\tenbit=cmssi10
\font\tenbbb=msbm10
\font\ninei=cmmi9 
\font\ninerm=cmr9

\font\ninesy=cmsy9 
\font\ninebbb=msbm10 at 9 pt
\font\eighti=cmmi8
\font\eightrm=cmr8

\font\eightsy=cmsy8
\font\eightbbb=msbm7 at 8 pt
\font\seveni=cmmi7 
\font\sevenrm=cmr7
\font\sevensy=cmsy7
\font\sevenbbb=msbm7
 
%
% large fonts
%

%

\newfam\bbbfam

% this macro defines a tenpoint font
\def\tenpoint{
\def\rm{\fam0\tenrm}
\textfont0=\tenrm \scriptfont0=\eightrm \scriptscriptfont0=\sevenrm
\textfont1=\teni \scriptfont1=\eighti \scriptscriptfont1=\seveni
\textfont2=\tensy \scriptfont2=\eightsy \scriptscriptfont2=\sevensy
\textfont3=\tenex \scriptfont3=\tenex \scriptscriptfont3=\tenex

\textfont\itfam=\tenit
\def\it{\fam\itfam\tenit}

\textfont\slfam=\tensl
\def\sl{\fam\slfam\tensl}

\textfont\bffam=\tenbf
\def\bf{\fam\bffam\tenbf}

\textfont\ttfam=\tentt
\def\tt{\fam\ttfam\tentt}

\def\sc{\tensc}

\def\bit{\tenbit}

\def\bbb{\fam\bbbfam\twelvebbb}
\textfont\bbbfam=\tenbbb
\scriptfont\bbbfam=\eightbbb
\scriptscriptfont\bbbfam=\sevenbbb
 
\normalbaselineskip=12pt

\setbox\strutbox=\hbox{\vrule height10pt depth4pt width0pt}%
\normalbaselines\rm}

% this macro defines a twelvepoint font

\def\twelvepoint{
\def\rm{\fam0\twelverm}
\textfont0=\twelverm \scriptfont0=\ninerm \scriptscriptfont0=\sevenrm
\textfont1=\twelvei \scriptfont1=\ninei \scriptscriptfont1=\seveni
\textfont2=\twelvesy \scriptfont2=\ninesy \scriptscriptfont2=\sevensy
\textfont3=\tenex \scriptfont3=\tenex \scriptscriptfont3=\tenex

\textfont\itfam=\twelveit
\def\it{\fam\itfam\twelveit}

\textfont\slfam=\twelvesl
\def\sl{\fam\slfam\twelvesl}

\textfont\bffam=\twelvebf
\def\bf{\fam\bffam\twelvebf}

\textfont\ttfam=\twelvett
\def\tt{\fam\ttfam\twelvett}

\def\sc{\twelvesc}
\def\bit{\twelvebit}

\def\bbb{\fam\bbbfam\twelvebbb}
\textfont\bbbfam=\twelvebbb
\scriptfont\bbbfam=\ninebbb 
\scriptscriptfont\bbbfam=\sevenbbb
 
\normalbaselineskip=14pt

\setbox\strutbox=\hbox{\vrule height10pt depth4pt width0pt}%
\normalbaselines\rm}

% this macro defines a fourteen point font

\def\fourteenpoint{
\def\rm{\fam0\fourteenrm}
\textfont0=\fourteenrm \scriptfont0=\twelverm \scriptscriptfont0=\tenrm
\textfont1=\fourteeni \scriptfont1=\twelvei \scriptscriptfont1=\teni
\textfont2=\fourteensy \scriptfont2=\twelvesy \scriptscriptfont2=\tensy
\textfont3=\tenex \scriptfont3=\tenex \scriptscriptfont3=\tenex

\textfont\itfam=\fourteenit
\def\it{\fam\itfam\fourteenit}

\textfont\slfam=\fourteensl
\def\sl{\fam\slfam\fourteensl}

\textfont\bffam=\fourteenbf
\def\bf{\fam\bffam\fourteenbf}

\textfont\ttfam=\fourteentt
\def\tt{\fam\ttfam\fourteentt}

\def\sc{\fourteensc}
\def\bit{\fourteenbit}

\def\bbb{\fam\bbbfam\twelvebbb}
\textfont\bbbfam=\fourteenbbb
\scriptfont\bbbfam=\tenrm 
\scriptscriptfont\bbbfam=\eightrm

\normalbaselineskip=16pt

\setbox\strutbox=\hbox{\vrule height10pt depth4pt width0pt}%
\normalbaselines\rm}

% Set some default dimensions
%
\twelvepoint

\abovedisplayskip 14pt plus 3pt minus 10pt%
\belowdisplayskip 14pt plus 3pt minus 10pt%
\abovedisplayshortskip 0pt plus 3pt%
\belowdisplayshortskip 8pt plus 3pt minus 5pt%
\parskip 3pt plus 1.5pt
\hsize=6.5in
\vsize=8.9in

\def\Bbb{\bbb}

%
% useful macros

\fi			% end of \ifx\UsualIsLoaded

\def\SBIMSMark#1#2#3{
 \font\SBF=cmss10 at 10 true pt
 \font\SBI=cmssi10 at 10 true pt
 \setbox0=\hbox{\SBF Stony Brook IMS Preprint \##1}
 \setbox2=\hbox to \wd0{\hfil \SBI #2}
 \setbox4=\hbox to \wd0{\hfil \SBI #3}
 \setbox6=\hbox to \wd0{\hss
             \vbox{\hsize=\wd0 \parskip=0pt \baselineskip=10 true pt
                   \copy0 \break%
                   \copy2 \break% 
                   \copy4 \break}}
 \dimen0=\ht6   \advance\dimen0 by \vsize \advance\dimen0 by 8 true pt
                \advance\dimen0 by -\pagetotal
 \dimen2=\hsize \advance\dimen2 by .25 true in
%
%   Check for publication info
%
%  \newread\jref
  \openin2=publishd.tex
  \ifeof2\setbox0=\hbox to 0pt{}
  \else 
     \setbox0=\hbox to 3.1 true in{
                \vbox to \ht6{\hsize=3 true in \parskip=0pt  \noindent  
                \input publishd.tex 
                \vfill}}
  \fi
  \closein2
  \ht0=0pt \dp0=0pt
 \ht6=0pt \dp6=0pt
 \setbox8=\vbox to \dimen0{\vfill \hbox to \dimen2{\copy0 \hss \copy6}}
 \ht8=0pt \dp8=0pt \wd8=0pt
 \copy8
 \message{*** Stony Brook IMS Preprint #1, #2 ***}
}

\font\de=cmssi12

%\magnification\magstep1

\def\gm{\gamma}
\def\gmn{\gamma_N}
\def\Gm#1#2{\gm_{|_{[#1, #2]}}}

\def\GM{\Gamma}

\def\td{\tilde}
\def\ro#1#2{\rho (\gm;#1, #2)} %or whatever else we want to call it.

\def\ron#1#2{\rho ({\gm_N};#1,#2)}

\def\Lg#1#2{\ell(\gm;{#1}, {#2})}
\def\Ac#1#2{A(\gm;{#1}, {#2})}
\edef\Dist#1#2{d(\gm ({#1}), \gm({#2}))}

\def\del{\partial}
                              %Action

\def\circle{S^1}
\def\ra{\rightarrow}
\def\reals{{\Bbb R}}    

\def\P{TM\times \circle}
\def\tP{T\td M\times\circle}
\def\Man{Ma\~n\'e}
\def\Mane{Ma\~n\'e}
\def\tM{\td M}
\def\Poin{{{\Bbb H}^n}}
\def\QED{  \rlap{$\sqcup$}$\sqcap$ }
\def\Geos{{\cal G}}
\def\Min{{\cal M}}
\def\bs{{\bar\sigma}}
\def\tp{{\tilde\phi}}
\def\ie{{\it i.e.}}
\def\eg{{\it eg.}}
\def\AM{Aubry-Mather}
\def\cf{{\it cf}}
\def\cite#1{{\bf #1}}
\def\xb{{\bf x}}
\def\xd{{\dot x}}
\def\T{{\bf T}}
\def\R{{\bf R}}
\def\dx{{\dot x}}
\def\pt{\phi_t}
\def\tz{{\tilde z}}
\def\tz{{\tilde z}}
\def\tpt{{\tilde\phi_t}}
\def\bibitem#1{\item{[#1]}}

\def\RR{{\Bbb R}}
\def\NN{{\Bbb N}}

\def\HH{{\Bbb H}}

   % macro with one variable [parameter] \abs x -> |x|
\def\norm#1{\left\Vert #1\right\Vert}   
% macro with one variable [parameter] \norm x -> ||x||
\medskip

\medskip
\def\Cross{\bigm| \kern-5.5pt \not \ \, }
\def\cross{\mid \kern-5.0pt \not \ \, }
\medskip

\medskip
\hyphenation{pa-ra-meters}
\hyphenation{pa-ra-meter}
\hyphenation{lem-ma}
\hyphenation{lem-mas}
\hyphenation{to-po-logy}
\hyphenation{to-po-logies}
\hyphenation{homo-logy}

\SBIMSMark{1996/1a}{January 1996}{}
\centerline{\bf Lagrangian Systems on Hyperbolic Manifolds}
\bigskip
\centerline{Philip Boyland and  Christophe Gol\'e}

\bigskip

{\narrower

\noindent{\bf Abstract:}
 This paper gives two results that show that the
dynamics of a time-periodic Lagrangian system on a hyperbolic manifold
are at least as complicated as the geodesic flow of  a hyperbolic
metric. Given a hyperbolic geodesic in the Poincar\'e ball, Theorem
A asserts that there are minimizers of the lift of the 
Lagrangian system that are a bounded distance away and have a  
variety of approximate speeds. Theorem B gives the existence of a
collection of compact invariant sets of the Euler-Lagrange flow that are
semiconjugate to the geodesic flow of a hyperbolic metric.
These results can be viewed as a generalization of
the \AM\ theory of twist maps and the Hedlund-Morse-Gromov theory of 
minimal geodesics on closed surfaces and hyperbolic manifolds.
\bigskip
}

\noindent{\bf Section 0: Introduction.}

The notion of stability in Dynamical Systems
refers to  dynamical behavior that persists
under perturbation. Stability under small perturbations is perhaps best known, but dynamical persistence under large
perturbations (in a restricted class) is often studied and has
proved to be quite powerful. Large perturbation theories 
usually have a strong topological component. This is  because  
behavior  that persists under  large 
perturbations must be very fundamental to the system, and 
the most fundamental aspect of a dynamical 
system is the topology of the underlying manifold.
In applying stability results, one usually begins with a
model system whose dynamics are understood and then 
perturbs it. The stability theorems indicate which dynamics of the 
model system must be present in the perturbed system.
This provides a framework for the investigation the other dynamics
present in the perturbed system.

This paper presents  stability results for the dynamics of time-periodic
Lagrangian (or Hamiltonian) systems for which the 
configuration manifold carries a hyperbolic metric, \ie\ a 
metric of  constant negative curvature. In this case the 
model system is the geodesic flow of a hyperbolic metric.
The results  generalize and/or are closely related to several  
theories that contain what may be viewed as stability results, 
for example,  the \AM\ theory of twist maps  and the  
Hedlund-Morse theory of minimal geodesics on closed surfaces.
The connection between geodesic flows, 
Euler-Lagrange flows, and the \AM\ theory has been  explored in 
 [B1], [BK], [BP], [Ma1],  [Mo] and elsewhere.
Our work also builds on a generalization of these theories
due to  Mather [Ma1], [Ma2] (see also [Mn]).
(For a complete survey of the connection of the results here to various other theories, see [BG].)

These theories share the property that the orbits of the 
dynamical system under consideration correspond to extremals of a 
variational problem defined in the universal cover of the 
configuration space. The orbits that correspond to minimums of the 
variational problem have special properties; they  behave 
approximately like the solutions to the  variational problem 
associated with the model system. This is natural because all orbits of the model problem are  minimizers. 
 
Theorem A gives our  first way of formalizing the idea that  
time-periodic Lagrangian systems on hyperbolic manifolds are at 
least as complicated as the geodesics of a hyperbolic metric.
Given a hyperbolic geodesic in the Poincar\'e ball, $\Poin$, the 
theorem asserts 
that there are minimizers of the Lagrangian system that are a 
bounded distance away and have a  variety of approximate speeds. 
Given a path
$\gm:\reals\ra\Poin$, the notation $\rho(\gm;a,b)$ means the 
average displacement in $\Poin$ over the time interval 
$[a,b]$., \ie\ the distance from $\gm(a)$ to $\gm(b)$ divided 
by $b-a$.

\medskip

{\bf Theorem A:} {\it Let $(M,g)$ be a closed hyperbolic 
manifold. 
Given a Lagrangian $L$ which satisfies Hypothesis 1.0,  there are 
sequences $k_i, \kappa_i, T_i$ in $\RR^+$  depending only on  
$L$,  with $k_i$ increasing to infinity, such that,  for any 
hyperbolic geodesic $\Gamma \subset \HH^n =\tilde M $, there are 
minimizers 
$\gm_i:\RR\to\td M$ with $d(\gm_i, \GM) \leq \kappa_i $, 
$\gm_i(\pm \infty)= \GM(\pm \infty)$,  and 
 $k_i\leq \rho(\gm; c,d)\leq k_{i+1} $ whenever $d-c\geq T_i$. 
}
\medskip

The basic idea of the proof is a limit argument that goes back to 
Morse [M]. Given a hyperbolic geodesic and a speed,
we approximates the geodesic by a long minimizing segment with 
the correct average speed. We then let the approximating segment 
get longer and longer and pass to the limit.  In order to pass to this limit we need some uniform control on the speed and geometry of 
the minimizing segments. This control comes from showing (Prop 2.1) that minimizing segments are quasi-geodesics in the sense of Gromov. Further, the quasi-geodesic constants depend only on the average speed of the minimizer.

Exact symplectic twist maps on the cotangent bundle of a manifold 
(defined in \S 1.3) are in many ways the discrete analogs of the E-L flow. Throughout the paper we indicate how the results for
Lagrangian systems can be adapted for the twist map case. In particular,
there is a twist map version of Theorem A.

The second main result, Theorem B, 
focuses on the dynamics generated by the Lagrangian
and is a kind of globalization of 
Theorem A. One way to formulate the fact that a perturbed system 
is at least as complicated as the model system (\ie\ the dynamics of
the model systems don't go away) is to show that the perturbed 
system always has a invariant set that carries the dynamics of  the 
minimal model. More precisely, one shows that there is a
compact invariant set that is semiconjugate to the minimal model. 
MacKay  and Denvir ([MD]) have recently extended Morse's results to the 
case with boundary and proved a result giving this semiconjugacy. 
Also Gromov ([G1]) and others have obtained semiconjugacies in the case of geodesic flows.

Theorem B is a semiconjugacy theorem for time-periodic 
Lagrangian systems. Given such a  Lagrangian, 
its Euler-Lagrange equations yield a second order time-periodic
differential equation on the tangent bundle $TM$,  and thus a vector 
field on $\P$. The solution flow of the vector field (when it exists) is called the E-L flow. 
The set  $\Min\subset\P$ consists of  all the orbits that correspond to 
minimizing paths in the  configuration space.

\medskip
{\bf Theorem B:} {\it Let $(M, g)$ be a closed hyperbolic manifold 
with geodesic flow $g_t$. Given a
Lagrangian $L$ which satisfies Hypotheses 1.0   with
E-L flow $\phi_t$, there exists  sequences $k_i$  and $T_i$  
with $k_i$ increasing to infinity, and a family  of compact, 
$\phi_t$-invariant  sets $X_i \subset \Min$  so that
for all $i$, $(X_i, \phi_t)$ is semiconjugate to $(T_1 M, g_t)$ 
and  $k_i\leq \rho(\phi_t(x);0,T)\leq k_{i+1}, $
whenever $T\geq T_i$ and $x\in X_i$.}
\medskip

Note that the geodesic flow of a hyperbolic metric is 
transitive Anosov and is  therefore  Bernoulli, 
has positive entropy, etc (see, \eg\ [HK]).
Thus Theorem B implies that the E-L flow is always dynamically very complicated.

In the last subsection we remark on how the semiconjugacy in Theorem B 
can be used
to find ergodic $\phi_t$-invariant measures for which an average speed exists almost everywhere. Further, each of these measures  ``shadow'' a unique ergodic
$g_t$-invariant measure.

\medskip
\noindent{\bf Section 1: Preliminaries.}
 
In this section we introduce notation and recall some basic results 
needed in the sequel. For a general discussion of  Lagrangian 
systems see [AM]. For  thorough discussions of Lagrangian 
systems and minimizers the reader is urged to consult 
[Ma1], [Ma2],   and [Mn]. For more details on symplectic twist maps, 
the reader is referred to [Gl] or [MMS]
(\cf\ [BK] and [K])

\medskip
\noindent{\bf\S\ 1.1 Lagrangian systems.}
The main objects in the 
Lagrangian formulation of mechanics are a {\de configuration manifold}
$M$ and  a real valued function 
called a {\de Lagrangian} defined on the tangent bundle $TM$.
The  configuration spaces of interest here are  closed manifolds $M$ 
with a fixed Riemannian metric $g$. The induced norm on the 
tangent bundle is denoted $\norm{v}$. We consider time-periodic 
systems determined by a  $C^2$-Lagrangian $L: TM\times \circle\to 
\RR$.  The basic variational problem is to  find curves $\gm:[a,b]\to 
M$ that are extremal for the {\de action}
$$
A(\gm)=\int_a^b L(\gm,\dot\gm,t) dt
$$
among all absolutely continuous curves $\beta:[a,b]\to  M$ 
that have the same endpoints 
$\beta(a)=\gm(a), \beta(b)= \gm(b)$.

Under appropriate hypothesis (eg. $\gm$ is $C^1$),
such a $\gm$ satisfies
the Euler-Lagrange second order differential equations
$$
{d\over dt}{\del L\over \del v}(\gm(t),\dot \gm(t),t)-{\del L\over 
\del x}(\gm(t),\dot
\gm(t),t)=0. 
$$
Using local coordinates these equations yield a 
first order time-periodic differential equation
on $TM$, and thus in the standard way, a vector field on
$\P$.  Since $\P$ is not compact it is possible that trajectories of this 
vector field are not defined for all time in $\RR$ and
thus do not fit together to give a  global flow (i.e. an $\RR$-action). 
When the  
flow does exist, it is called the {\de Euler-Lagrange (E-L) flow}.

All Lagrangians in this paper are assumed to 
satisfy the following hypotheses. 

\medskip
\noindent{\bf  Standing Hypotheses 1.0:}
{\it $L$ is a $C^2$-function $L:\P\ra\RR$ that satisfies:
 
(a) {\it Convexity:} $\del^2 L\over \del v^2$ is positive definite.
 
(b) {\it Completeness:} The Euler-Lagrange flow determined by L 
exists.
 
(c) {\it Superquadratic:} There exists a $C > 0$ so that $L(x,v,t)\geq 
C\norm{v}^2$.
}

\medskip
These assumption are the same as those in Mather [Ma1]
and \Mane\ [Mn]
apart from (c), where they only assume
${L(x,v,t)\over \norm v} \to\infty$
when $\norm{v}\to +\infty$. We need the stronger condition  in  
the proof of Lemma 2.2.
Note that the addition of a constant to $L$ 
does not change the E-L flow. 

{\bf Remark:} Condition 1.0(a) is the classical Legendre
condition. Thus the E-L flow derived from such an $L$ is conjugate under
the Legendre transform to a Hamiltonian flow on $T^*M\times S^1$.
In particular, if $L$ is time {\it independent}, orbits of the E-L flow are 
constrained to the Legendre transforms of the constant energy surfaces of the Hamiltonian flow. 
In the time independent case, the use of this fact along with the
Jacobi metric (or finsler) results in a considerable simplification of many the proofs in this paper. The time dependent case is more much interesting
and is our focus here.

\medskip
{\bf Example: Mechanical Lagrangians.} 
As pointed out by \Mane, Hypotheses 1.0  are satisfied for 
mechanical Lagrangians, i.e. those of the  form 
$$
L(x,v,t)= {1\over 2}\norm{ v}_h^2 -V(x,t),\ \  V\leq 0
$$
where the norm is taken with respect to any Riemannian metric $h$.
(In fact,  under some conditions, one may allow the norm to vary with time. See 
[Mn], page 44). 
\medskip

\noindent{\bf \S 1.2 Minimizers.}
Of particular interest in  the Lagrangian theory
 are extremals of the variational problem  
that minimize in the following sense.  
If $\td M$ is a regular covering space of $M$,
 $L$ lifts to a real valued function (also called $L$)  defined on $\tP$.
A curve segment $\gm:[a,b]\to \tM$ is called a $\tM$-{\de minimizing 
segment} or an $\tM$-{\de minimizer} 
if it minimizes
the action among all absolutely continuous curves $\beta:[a,b]\to 
\tM$
which have the same endpoints.

A fundamental theorem of Tonelli implies that
 if $L$ satisfies Hypotheses 1.0,
then given  $a < b$ and two distinct points
$x_a,x_b\in\tM$ there is always a minimizer $\gm$ with  $\gm(a) = 
x_a$
and $\gm(b) = x_b$. Moreover such a $\gm$ is automatically $C^2$ and
satisfies the Euler-Lagrange equations  (this uses the completeness of
the E-L flow). Hence its lift, $ (\gm(t),\dot\gm(t), t)$,
to $\P$ is a solution of the E-L flow. 
A curve  $\gm:\RR\ra\tM$ is called a {\de 
minimizer} if $\Gm ab$ is a minimizer for all $[a,b]\subset\RR$. 
When the domain of definition of  a curve is not explicitly given, it is 
assumed to be $\RR$. 

Mather [Ma1] and \Man\   [Mn]
use $\overline M$ minimizers where
 $\overline M$ is the universal free Abelian cover. The universal cover
(which we denote $\tM$ from now on) is used here.  If $\gm$
 is an $\tM$-minimizer, we will simply say it is a minimizer.

Our main task is to get control of the speed and geometry of 
minimizers. 
Given a smooth curve $\gm:[c,d]\ra\tM$  and a segment $[a,b]\subset 
[c,d]$,  the {\de average
displacement} in the  cover over the interval $[a,b]$ is measured by 
$$
\ro  ab= {d (\gm(a), \gm(b))\over b-a}
$$
where $d$ is the topological metric on the universal cover 
constructed from the lift of the given Riemannian 
metric $g$. The length of
$\Gm ab$  is denoted $\Lg ab$, and the {\de action}
over the interval $[a,b]$ is
$$
\Ac ab = \int_a^b L(\gm,\dot\gm,t) dt.
$$
In all these notations  the absence of the last two arguments 
indicates
the quantity is computed for the entire interval of definition, thus 
$\rho(\gm) = \ro cd $, etc.

Using the fact that $L$ is superquadratic as assumed in 1.0(c), 
we obtain 
simple but very useful estimates on the average action of minimizers. 
The estimates are essentially in \Mane\ [Mn]  and 
Mather [Ma1], but 
the versions given here are  slightly more exact as  
our assumptions on $L$ are slightly stronger.
 The proof follows  [Mn], Theorem 3.3. 

\medskip

{\bf Lemma 1.1:} {\it  
Given a Lagrangian $L$ satisfying  Hypothesis 1.0, let 
$C_K^{max}=$\hfill\break ${1\over K} \sup\{ L(x,v,t) : \norm{v}\leq K\}$ and 
$C_K^{min} = {C K \over 4}$, where $C$ is the constant in 1.0(c).  If 
$\gm$ is a minimizer and 
$ \ro ab =K$, then $$C_K^{min}K\leq {1\over b-a}{\Ac ab}  \leq 
C_K^{max}K. 
$$
}

\medskip {\bf Proof:}  
 If $\GM:[a,b]\to \td M$ with 
$\gm(a)=\GM (a)$ and $\gm (b)=\GM (b)$ is a minimizing {\it 
geodesic} segment with respect to the given metric $g$, then 
$\norm{\dot \GM}= \rho(\GM;a,b)=\ro ab 
=K$. Thus, 
$$
A (\gm)\leq A(\GM)\leq\int_a^b K C_K^{max}\; dt= 
K C_K^{max} (b-a),
$$
yielding the upper bound.

For the lower bound, first note that
$$\eqalign{
K(b-a)=\Dist ab 
& \leq \int_a^b \norm{\dot \gm}\;dt
=\int_{\norm{\dot \gm}>{K\over 2}}\norm{\dot \gm}\;dt 
+\int_{\norm{\dot \gm}
\leq {K\over 2}}\norm{\dot \gm}\;dt\cr 
&\leq \int_{\norm{\dot \gm}>{K\over 2}}\norm{\dot \gm}\;dt + 
{K\over 2}(b-a),
\cr}
$$

and so
$$\int_{\norm{\dot \gm}>{K\over 2}}\norm{\dot \gm}dt\geq 
{K\over 2}(b-a). 
$$  
Thus using the Cauchy-Schwartz inequality and 1.0(c),
$$\eqalign{
A(\gm)
&\geq \int_{\norm{\dot \gm}>{K\over 2}}L(\gm,\dot \gm,t)\;dt
\geq C\int_{\norm{\dot \gm}>{K\over 2}}\norm{\dot \gm}^2\;dt\cr
&\geq {C\over b-a}\left( \int_{\norm{\dot \gm}>{K\over 
2}}\norm{\dot \gm}\;dt\right)^2
\geq {CK^2\over 4}(b-a).\cr}
$$
\hbox{} \QED

\medskip
{\bf Remark: } Note that  $K\mapsto C_K^{min}$ is a continuous  
function  that increases monotonically to infinity, while $K\mapsto 
C_K^{max}$ is  continuous and grows to infinity (since $C_K^{max} 
\geq C K$).  Note also that $K 
C_K^{max}$ is monotone nondecreasing. These facts will be 
used frequently in the sequel without further mention.
\medskip

\medskip
{\bf Example: Mechanical Lagrangians.} 
Consider again the  mechanical Lagrangian
$$
L(x,v,t)= {1\over 2}\norm{v}_h^2 -V(x,t)$$ where
$\norm{\ }_h$ comes from a Riemannian metric, and we
let $\max V =0$ and $\min V := V_{min}$.
If $B_1$ and $B_2$ are the positive constants such that
$$
B_1\norm{v}^2\leq \norm{v}_h^2\leq B_2 \norm{v}^2,\eqno{(1.1)}
$$ 
where $\norm{\cdot}$ is the norm coming from
the fixed reference metric $g$ (\eg\ of constant curvature), one readily computes that 
$$
C_K^{min}={B_1K\over 8}, \quad C_K^{max}= 
{1\over 2} B_2 K-{V_{min}\over K}
$$
(see [BG] for slightly better estimates in mechanical case).
\medskip

{\bf \S 1.3 Exact symplectic twist maps.}
An {\de exact symplectic twist map} $F$  is a map from
a subset $U$ of the cotangent bundle  of a manifold $N$ (which we
allow to be noncompact) into $U$,
which comes equipped with a {\de generating function} $S:N \times N \ra \reals$ that satisfies
$$
F^*(p\;dx)-p\;dx= P\;dX-p\;dx= dS(x, X),\eqno{(1.2)}
$$
where $(X,P)$ are the coordinates of $F(x, p)$ (this can also be
written in a coordinate free manner). 

Because the one-form
$P\;dX-p\;dx$ in (1.2) is exact, one says that $F$ is {\de exact}.  
Note that taking the exterior differential of (1.2) yields
 $dP\wedge dX= dp\wedge dx$, and so any exact $F$ is also symplectic,
\ie\ it preserves the standard symplectic form.
The fact that  $S$ is expressed using  the coordinates  $(x,X)$  instead of $(x,p)$ is the {\de twist condition}.  
Given $S$, one can retrieve the map (at least implicitly)
from $p=-{\del S\over \del x} and P={\del S\over \del X}$.
This can be done globally (\ie\   $U=T^*N$) only
when $N$ is diffeomorphic to a fiber of $T^*N$, for example  when
$N$ is the covering space of the n-torus or of a manifold of constant negative
curvature.

The variational problem for Lagrangian systems translates into a  discrete variational problem for twist maps: the role of  curves in the continuous setting is taken by sequences of points (``integer time curves''), and the
action of a finite sequence $\xb=\{x_n, \ldots, x_m\}$
is given by $W(\xb)= \sum_n^{m-1} S(x_k,x_{k+1})$.
This corresponds closely to the continuous setting 
when the exact symplectic twist map $F$ is the time-one map of
an E-L flow. In this case,  $S(x,X)=\int_0^1L(x, \xd,t)dt$,
where $x(t)$ is the minimizer with endpoints $x$ and $X$.

In direct correspondence to  Lagrangian systems, critical points of $W$ 
(with fixed time and configuration endpoints) correspond to  orbits of $F$
(this is closely related to 
the method of  broken geodesics  in Riemannian geometry).
Action minimizers are sequences that minimize $W$ over any of their
subsegments.  The natural growth condition on the generating function 
$$
S(x, X)\geq C \; dist^2 (x, X), 
$$  
implies the  analog of Tonelli's theorem: minimizers always exist 
between any two points over any given (integer) interval of time. 
Moreover, there is 
an exact analog of Lemma 1.1: the average action of minimizers
is bounded below and above by functions of the average
displacement. The proof is virtually identical to the continuous time case, 
replacing geodesics with orbits of the time-one of the geodesic flow.

\medskip
{\bf Example: Generalized standard maps.}
 Let M be $\T^n$ or a closed hyperbolic manifold, and
let $N = \tM$ be the universal cover $\R^n$ or $\Poin$,
respectively. 
On the covering space, define the  
{\de generalized standard map} using its generating function  
$\tM\times\tM\ra\reals$,
$$S(x, X)={1\over 2} dist^2(x, X) +V(x)$$
where the distance $dist$ is  induced by the Euclidean metric in $\R^n$, 
or the hyperbolic metric on $\Poin$, and $V(x)$ is $\pi_1(M)$-equivariant,
\ie\  it descends to a function on $M$. A short argument shows that
one can use the relation (1.2) to solve for $(X, P)$ in terms of
$(x,p)$ and thus obtain
an exact symplectic twist
map on $T^*\tM$ that, in turn, 
 induces a map on $T^*M$ (also called a twist map).
For more general examples, \cf\ [Gl].
 
{\bf Remark:} In certain cases the twist map theory overlaps with the
continuous theory. If a twist map $f$ of  $T^*\T^n$ has  a 
 generating function  that is
super quadratic in $\norm{X-x}$, the mixed  partial 
$ \partial_{12} S$ is symmetric,
and for some $a>0$ satisfies the convexity condition
$$
< \partial_{12} S (x, X).v, v >\  \leq -a \norm{v}^2
$$
uniformly in $(x, X)$, then $F$ is the  time-one map of an
E-L flow derived from a 
one-periodic Lagrangian that is superquadratic in the velocity. 
Moser [Mo] gives the proof in the case $n=1$. Bialy and
Polterovitch remark in  [BP]
that Moser's proof goes through in the case $n>1$.  This 
is not quite so, but they subsequently obtained a different
proof (personal communication).  Note that the generating function
for the generalized standard map satisfies these hypothesis.

\medskip

\noindent{\bf \S1.4 Hyperbolic Geometry.}
We recall some basic facts about  hyperbolic geometry and 
manifolds. For more information see \eg\ [BKS].
 A closed manifold $M$ is called {\de hyperbolic} if there 
is a Riemannian metric $g$ on $M$ that has  curvature
 identically equal to $-1$. The universal cover $\tM$ of a $n$-dimensional
 hyperbolic manifold is homeomorphic to $\RR^n$. We identify
$\tM$ with the n-dimensional Poincar\'e Disk $\Poin$, and so the 
group of covering transformations  can  be identified with a discrete 
subgroup, isomorphic to $\pi_1(M)$,
of the set of isometrics of $\Poin$. This group action has a 
fundamental domain with compact closure and  under 
the quotient by the action  the metric on $\Poin$  descends
to the hyperbolic metric on $M$. As a consequence of the Mostow Rigidity Theorem, a closed hyperbolic manifold of dimension three or greater carries a unique hyperbolic metric.  If $M$ is a surface, it carries  many hyperbolic metrics, but their geodesic flows are all orbit equivalent ([G1])

The sphere at infinity, $S_\infty$, is the usual Euclidean
sphere that is the boundary of $\Poin$ in the {\it Euclidean} 
topology.
The geodesics in $\Poin$ are semi-circles that are orthogonal
to $S_\infty$. In this paper these hyperbolic geodesics will always
be oriented and parameterized by arclength. For a geodesic
$\Gamma$, the notations $\Gamma(\infty)$ and $\Gamma(-\infty)$
refer to the  limit points of 
$\Gamma$ on $S_\infty$ in the forward and backward directions,
respectively. More generally, for $\gm:\RR\ra\Poin$, 
the notations $\gm(\infty)$  and $\gm(-\infty)$ refer to the limit 
points {\it in the Euclidean topology} of
$\gamma$ in forward and backward time if these limits
exist and are contained in $S_\infty$. Implicit in this notation is the
fact that $\gm$ has no nontrivial limit points in $\Poin$,
 {\it i.e.} $\gm:\RR\ra\Poin$ is proper.

For a pair of points $x,y\in\Poin$, $d(x,y)$ denotes their
distance in the hyperbolic metric. The notion of a  quasi-geodesic,  central to Gromov 's work on  hyperbolic groups  is also of central importance here. Given 
$\lambda > 1$ and $\epsilon>0$,  a curve $\gamma:\RR\ra
\Poin$ or a curve segment $\gamma:[a,b]\ra\Poin$ is called
a $(\lambda,\epsilon)$-quasi-geodesic if
$$
\lambda^{-1}(d-c) - \epsilon \leq d(\gm(c),\gm(d)) \leq
\lambda (d-c) + \epsilon
$$
 for all $[c,d]$ in the domain of $\gm$. 

The next theorem, usually called ``Stability of
quasi-geodesics'',  gives the most important property of
quasi-geodesics. It is true in the broader  context of 
what are usually called $\delta$-hyperbolic spaces, but we just state the 
result needed here. Given two closed subsets $X,Y\subset \Poin$,
$d(X,Y)$ denotes their Hausdorff distance  as induced by the 
hyperbolic
metric. For a proof and more information see [GH], [CDP] or [G2].

\medskip

{\bf Theorem 1.2:} {\it Given  $\lambda > 1$ and $\epsilon > 0$, 
there
exists a $\kappa>0$ so that whenever $\gm$ is a $(\lambda, 
\epsilon)$-quasi-geodesic segment  in $\Poin$ and $\Gamma_0$ is 
the geodesic 
segment connecting
the endpoints of $\gm$, then $d(\gamma, \Gamma_0) <\kappa$.
If $\gm$ is a $(\lambda, \epsilon)$-quasi-geodesic, then 
$\gm(\infty)$ and $\gm(-\infty)$ exist and further, if
$\Gamma$ is the geodesic connecting
$\gm(\infty)$ and $\gm(-\infty)$, then 
 $d(\gamma, \Gamma) < \kappa$.
}

\medskip
{\bf Remark:}  Although the  notion of quasi-geodesic
 makes sense on any Riemannian manifold, or even metric space,  it only
yields the strong consequence as in the last theorem for hyperbolic manifolds 
or, more generally,  $\delta$-hyperbolic spaces. 
One can easily construct counter examples to the theorem in  Euclidean space. 
These counter-examples contain the seeds of failure for the analog of 
Theorem A on the 3-torus. See  Section 3.3 in [BG] for more details.

\medskip

\noindent{\bf Section 2: Minimizers and  quasi-geodesics.}

The main purpose of this section is to prove Theorem A. 
\medskip
\noindent{\bf \S 2.1 Minimizing segments are quasi-geodesics.}
Throughout this section we fix a Lagrangian  $L$ that satisfies the 
Hypotheses 1.0. The first proposition  gives
uniform upper and lower bounds on the local 
average displacement of minimizing 
segments with a given total average displacement in the cover.  Part (a),
due to Mather (see the proof of Proposition 4 in [Ma1]),
gives an absolute upper bound of the velocity.  Mather considers 
minimizers in the universal  free Abelian cover, but  his proof 
works without change in the universal   cover.
Part (b)  says that points on minimizers 
cannot go too slow for too long. The main ingredient in the proof of 
(b) is Lemma 2.2. Its proof was inspired by Mather's proof of 
(a).  It uses an argument of 
curve shortening type. One assumes that a minimizer 
does not have the  desired 
property, and this allows one can to construct another curve that has 
lesser action, yielding a contradiction. Part (c) is a consequence of (a) and (b).

\medskip
{\bf Proposition 2.1:} 
{\it Let $\gm:[a,b] \to \td M$ be a minimizing segment with average displacement $\ro ab  
=K $ and $b-a\in\NN$.

(a)  (Mather) There is a $K''>K$ depending only on $L$ and $K$ such 
that
 for all $t\in [a,b]$, $\norm{\dot\gm(t)}\leq K''$.

(b) There exists $K_0>0$   such that for all $K>K_0$ there is a $k''>0$  
and an $N_0\in \NN$ depending only on $K$ and  $L$ so that for any 
interval $[c,d]\subset  [a,b]$ 
with $d-c \ge N_0$, one has $\ro cd \geq k''$.

(c) With $K$ and $K_0$ as in (b), 
there are constants $\lambda > 1$ and $\epsilon>0$ 
depending only on $K$ and  $L$  so that 
 $\gm$ is a $(\lambda, \epsilon)$-quasi-geodesic segment.}
\medskip

{\bf Remark:} 
 This proposition is true on any compact Riemannian 
manifold for Lagrangian systems satisfying Hypothesis 1.0.
{\it If the manifold is hyperbolic}, Theorem 1.2
implies that a minimizer $\gamma$ stays at a uniform distance from a 
geodesic $\Gamma$. However, because Theorem 1.2 is  generally not true on non-hyperbolic manifolds, one cannot obtain a version of Theorem A
in that case. As noted after Theorem 1.2, 
this provides a heuristic explanation 
for  why  straightforward  generalizations of the Aubry-Mather theory 
fail on the 3-torus. In the Hedlund metric,
one can construct a sequence of minimizing segments, each of which are quasi-geodesics by Proposition 2.1, but whose distance to any geodesic grows to infinity (see [BG], Section 4.2).

\medskip

\noindent{\bf \S 2.2 The main technical lemma.} 
The main step in proving Proposition 2.1
is a technical lemma that deals with  a special case of
Proposition 2.1(b). It gives a lower bound on the average
displacement when the subinterval 
has a specified integer length.

\medskip
{\bf Lemma 2.2:} {\it  
There exists $K_0>0$   such that for all $K>K_0$ there is a $K'$ with 
$0<K' < K$ and an $N_0\in \NN$ so that whenever $\gm:[a,b] \to \td M$ 
is a minimizing segment with $b-a\in\NN$ and $\ro ab  =K $, then 
for  any interval $[c,d]\subset [a,b]$ with $d-c= N_0$, and $b-d\in 
\NN$,  
one has $ \ro cd \geq K' $.
}
\medskip
{\bf Proof:} 
Since $C_K^{max}$ is a continuous function of $K\geq 0$ that is bounded 
below by ${C K\over 4}$, $m := \min_{K\geq 0} C_K^{\max}$ is achieved at some finite 
$K$. If we let $K_0 = 28 m / C$ (with $C$ as in 1.0(c)), 
then for any $K > K_0$ we can find a $K'$ with $0 < K' <K$
and 
$$
0<7C_{2K'}^{max} < CK/4 = C_K^{min}.\eqno{(2.1)} 
$$ 
Now let $N_0$ be the positive even integer with 
$$
{K\over K'}\leq {N_0\over 2} <  {K\over K'}+1.\eqno{(2.2)}
$$
Given $K > K_0$ and $K'$ and $N_0$ as in the last paragraph, assume   
$\gm, a, b, c, d$ satisfy the hypothesis of the lemma, but $\ro cd <K'$. 
We will construct a curve $\gm^*:[a,b]\to \td M$ of lesser action than 
$\gm$, yielding a contradiction.

First note that we  can find $[{a'},{b'}]\subset [a,b]$ such that ${b'}-
{a'}=1,  \ {a'}-a\in \NN$, and $\ro {a'}{b'} >K$.
Indeed, since $\ro cd <K'$,   either $\ro ac >K$ or $\ro db >K$
(the lost speed must be made up for somewhere). Say $\ro ac >K$, the 
other case is similar.
Subdivide  $[a,c]$ into intervals of length 1 (recalling that  $b-a\in 
\NN$). Clearly, one of these
intervals, which  we call $[{a'},{b'}]$, will satisfy the above conditions. 
Note in particular that $d (\gm({a'}), \gm ({b'})) >K$.

We now construct $\gm^*:[a,b]\to \td M$.
Let ${b''}\in [{a'},{b'}]$ be such that $d( \gm({a'}), \gm({b''}))=K$, and let $n={N_0\over 2}$.
The curve $\gm^*$ is defined on  sub-intervals of $[a,b]$ as follows:
\medskip
\vbox{\settabs\+On: &[{b''}+n,c+n] aa, & :$\gm^*$ is the minimizing 
segment with $\gm^*({a'}) = \gamma({a'})$ and $\gm^*({b'' + n})= 
\gm^*({b''})$\cr
\+On: &$[a, {a'}]$, &  $\gm^*(t)=\gm(t)$\cr
\+&$[{a'}, {b''}+n]$,&  $\gm^*$ is a minimizing segment with 
$\gm^*({a'}) = \gamma({a'})$, $\gm^*({b'' + n})= \gm^*({b''})$\cr
\+&$[{b''}+n, c+n]$, &  $\gm^*(t)=\gm (t-n)$\cr
\+&$[c+n,d]$,& $\gm^*$ is a minimizing segment with $\gm^*(c+n)= 
\gm(c)$, $\gm^*(d) = \gm(d)$\cr
\+&$[d,b]$,&  $\gm^*(t)=\gm(t).$\cr}
\medskip
\noindent From the time periodicity of $L$, it follows that 
\medskip
\vbox{\settabs\+$A(\gm)-A(\gm^*)$&$\ =\ $& $A(\gm; {a'},{b''})$ & $+$
 &$A(\gm; c, d)$ & $-$
& $A(\gm^*;a', b''+n)$ & $-$ & $A(\gm^*;c+n, d)$ \cr
\+$A(\gm)-A(\gm^*)$&$\ =\ $& $A(\gm; {a'},{b''})$ & $+$
 &$A(\gm; c, d)$ & $-$
& $A(\gm^*;a', b''+n)$ & $-$ & $A(\gm^*;c+n, d)$ \cr
\+ &$\buildrel \rm def \over=$& $\ \ A_1$& $+$&$\ \ A_2$&
$ -$&$\ \ A^*_1$&$ -$&$\ \ A^*_2$.\cr}
\medskip
\noindent We will show that this difference is positive by deriving
 estimates for $A_1, A_2, A^*_1, $ and $ A^*_2 $.

To estimate 
$A_1$, note that  
$$
{d( \gm(a'),\gm(b'') \over b''-a'}={K\over b''-a'}\buildrel \rm 
def\over =\td K\geq K.
$$
Thus,  using Lemma 1.1,
$$
A_1= A(\gm;{a'},{b''}) \geq C_{\td K}^{min} \td K 
(b''-a')\geq C_K^{min} K.
$$
All we need about 
$A_2$ is that $A_2\geq 0.$

Now for $A^*_1$, note that 
$${d(\gm^*({a'}), \gm^*({b''+n}))\over b''+n -a'}  = {d( \gm({a'}), 
\gm({b''}))\over n+b''-a'}={K\over n+b''-a'}
<{K\over n}\leq K',$$
and hence using Lemma 1.1 and (2.2), 
$$
A^*_1= A(\gm^*;{a'},{b''+n})\leq C_{K'}^{max}K' 
(b''+n-a')\leq 3 K C_{2 {K'}}^{max}.
$$

Finally, to estimate $A^*_2$, we observe that 
$${d( \gm({c+n}), \gm({d}))\over d-(c+n)}= {d(\gm(c), 
\gm(d))\over n} ={d( \gm(c), \gm(d))\over d-c}\cdot {d-c\over 
n}<2K',
$$
and so using (2.2),
$$
A^*_2=A(\gm^*;{c+n},d)  \leq C_{2K'}^{max} 2K' (d-
c-n)= 2C_{2K'}^{max} K'n \leq 4C_{2K'}^{max} K.
$$
 Now, if $K$ and $K'$ satisfy (2.1), clearly
$$A(\gm)-A(\gm^*)=A_1+A_2-A^*_1-A^*_2\geq K(C_K^{min}-
7C_{2K'}^{max})>0,$$ 
 a contradiction to the assumption that $\gm$ is a minimizer.  
\QED\medskip

\noindent{\bf \S2.3 Proof of Proposition 2.1.} 
As noted above, part (a) is proved in [Ma1]. For part (b), first note 
that  Lemma 1.2 implies that whenever $d-c=N_0$ and $b-d\in \NN$, 
we have 
$$
\Lg cd\geq d(\gm(c), \gm(d))\geq K'(d-c).\eqno{(2.3)}
$$
We now show that  slightly weaker is true for more general $[c,d]$, 
specifically, 
that for   any $[c,d]\subset [a,b]$ with  $d-c>N_0$, 
$$
\Lg cd\geq {K'\over 2 }(d-c). \eqno{(2.4)}
$$

 Let $c'$ be the smallest integer translate of $a$ that is bigger than 
$c$ ( i.e. $c'= a+\lfloor c-a\rfloor +1$).
 Divide $[c',d]$ into $m$ intervals $[c_i, 
c_{i+1}]$ of length $N_0$, plus possibly
one interval of lesser length $r$ (i.e. $d-c'=mN_0 +r, r<N_0$). If we
let $d'= 
c'+mN_0= d-r$,
then
$$
{\Lg cd\over d-c}
%={\Lg c{c'} +\Lg{c'}{d'} +\Lg {d'}d\over d-c}
\geq {\Lg{c'}{d'}\over d-c}
$$
and 
$$
{\Lg{c'}{d'}\over d-c}={1\over d-c}\sum_0^{m-1}\Lg {c_i}{c_{i+1}}
 \geq {N_0\over d-c}\sum_0^{m-1} {\Dist {c_i}{c_{i+1}}\over N_0} 
\geq {mN_0\over d-c} K',$$
the last inequality coming from (2.3). But
$$
{mN_0\over d-c} K'= {mN_0\over mN_0 +r +(c-c')}K'>{mN_0\over 
(m+1)N_0 +1}K'\geq {mN_0\over (m+1)N_0}K'.
$$
Since ${m\over m+1}$ increases with $m$ and $m\geq 1$, we have 
$$
{\Lg cd\over d-c}\geq {1\over 2}K',$$
which proves (2.4).

 Now  using the Cauchy-Schwartz inequality and what we have just 
proved,
$$ \eqalign{
{1\over d-c} \int_c^d L(\gm,\dot\gm,t)dt
& \geq {C\over d-c} \int_c^d\norm{\dot \gm}^2dt 
\geq {C\over (d-c)^2} \left(\int_c^d\norm {\dot \gm}dt\right)^2\cr
&= {C\over (d-c)^2}\left(\Lg cd \right)^2 \geq {C(K')^2\over 4} \cr}
$$
On the other hand, part (a) shows that $\ro ab =K$, implies 
$\ro cd\leq K''$ for any $[c,d]\subset [a,b]$ and so by Lemma 1.1, 
$$
{1\over d-c} \int_c^dL(\gm,\dot\gm,t)dt\leq C_{K''}^{max}\ro cd 
={C_{K''}^{max}\Dist cd\over d-c}.
$$
Part (b) then follows by letting 
$$
k'' = {C(K')^2\over 4 C_{K''}^{max}}.
$$

Part (c) follows from (a) and (b) by letting
$\lambda = \max\{ K'', 1/k'', 1\}$ and $\epsilon = N_0/\lambda$. 
\QED

\medskip
\noindent{\bf \S2.4 Proof of Theorem A.}
Given a path $\gm:\RR\ra\tM$,  we denote by $D\gm= (\gm,\dot \gm)$ its differential, which is 
a path in $T\tM$. 
%We call $p:\tP\ra\Poin$ the projection from $\tP$ to $\Poin$ .

Fix an  oriented geodesic with a given parameterization by arclength
 $\GM:\RR\ra\Poin$ and a $K > K_0$ with $K_0$ as in Proposition 2.1.
Let $\gmn:[-N, N]\ra\Poin$ be a minimizing segment with
$\gmn(-N) = \GM(-K N)$ and  $\gmn(N) = \GM(K N)$, and
thus $\ron{-N}{N}=K$. By Proposition 2.1(c), $\gmn$ is a
$(\lambda,\epsilon)$-quasi-geodesic and so by Theorem 1.1,
\hfill\break
\noindent $d(\gmn,\GM([-KN,KN]))<\kappa$, with $\lambda, \epsilon$ and 
$\kappa$ 
depending on $K$ and not on $N$.

Fix a fundamental domain of the manifold $M$ in the
universal cover $\Poin$ and call it $S$.  We may assume that
$S$ has compact closure and that $\GM$ intersects $S$. If 
$S_\kappa$ denotes the closure of $\{ x: d(x,S) \leq \kappa \}$, 
then $\gmn$ intersects 
$S_\kappa$ for all $N$.
 Thus we may find a $t_N$ with $\gmn(t_N)\in S_\kappa$.

Now let $z_N = D\gmn(t_N)$ and let $\tau_N\in\circle$ be the residue of 
$t_N$ mod 
$1$. Then using the $K'' = K''(K)$ from
Proposition 2.1(a), we have that for all $N$, $(z_N,\tau_N)$ is
contained in the compact space 
$\{ (x,v,t)\in\tP: x \in S_\kappa \ \hbox{and}\ \norm{v}\leq K''\}$.
Thus there exits $(z, \tau)$ and a subsequence $N_i$, with 
$(z_{N_i},\tau_{N_i})\ra (z,\tau)$ in $\tP$. 

Let  $y = (z,\tau)$, and let $y(t)= (\gm(t), \dot \gm(t), t)$  be the trajectory of the 
lift of the E-L flow with initial conditions $y$.  If $k'', N_0, \lambda$ and
$\epsilon$  are as in Proposition 2.1,  then  the continuity of solutions 
of differential equations with respect to initial conditions implies that 
$\gm(t)$ 
satisfies $k'' \leq\ro cd\leq K''$  whenever $d-c \ge N_0$,
since each $\gamma_N$ satisfies this inequality. 
Further, $\gm$ is a $ (\lambda, \epsilon)$-quasi-geodesic. 
Thus by Theorem 1.1, there is an oriented geodesic 
$\Gamma_1$ with $d(\gm,\GM_1)<\kappa$ and 
$\gm(\pm\infty)=\GM_1(\pm\infty)$.  Finally, it is clear that 
$\GM_1 = \GM$ since the $\gmn$ converge to $\gm$ pointwise and 
all the $\gmn$ and $\gm$ are quasi-geodesics with the same 
constants. 

Doing this construction for each $K$  yields a family of minimizers. This 
family clearly contains the sequence  required for the theorem. \QED
\medskip

{\bf Remarks:}

 Note that in contrast to the \AM\ theory, Theorem A 
 does not yield all
speeds in every direction (in the \AM\ case there is just one 
direction). This is almost certainly not just artifact of the proof as there are autonomous mechanical systems on the two-torus
that have gaps in the speed spectrum (see Section 6 in [BG]). 
Nonetheless, one would expect that Theorem A could be improved.

 If $\GM$ is a closed geodesic, then the shadowing minimizer 
does not necessarily have to be a closed orbit. By minimizing the
action in the space of loops of integer period
in  each free homotopy class, one find periodic orbits of 
all free homotopy types (see [Bn] and Proposition 7.1 in [Mn])
 However, these periodic orbits may not be minimizers in the sense used here.

 The proof of analog of Theorem A for symplectic twist maps is 
virtually identical to that of the continuous case. By using the integers to
define a discrete
time, one has an analogous notion of 
quasi-geodesic for which the analog of Theorem 1.2 holds. Thus one obtains the Theorem using the twist map version of Proposition 2.1.

\medskip
\noindent{\bf Section 3: Semiconjugacy with the geodesic flow}

The previous section was concerned with minimizing 
curves in the universal cover and their relation 
to hyperbolic geodesics. In this
section we consider the dynamical implications of those results and 
focus on the E-L flow on $\tP$. We prove Theorem B 
by uniformizing the results of the previous
section and obtain semiconjugacies from subsets of the set of minimizers  
to the geodesic flow of the fixed hyperbolic metric.  

\medskip
\noindent{\bf \S 3.1 Definitions.}
The E-L flow $\phi_t$ on $\P$ lifts to a flow $\tilde\phi_t$ on $\tP$. 
There are two projections
that are of importance here. The projection from the covering space 
to
the base is $\pi:\tP\ra \P$. The projections from bundles
to the bases are both denoted $p$ and are $p:\tP\ra \tM$ and
$p:\P\ra M$.

Take an orbit $\tilde \phi_t(z)$ of the lift of the E-L flow and
let $\gm(t) = p(\tilde \phi_t(z))$. The notions that we used in 
 the previous section to describe curves $\gm$ will also be applied to 
orbits. Thus if $\gm$ is 
a minimizer or  a $(\lambda,\epsilon)$-quasi-geodesic, then this 
same
label is attached to the orbit. The subset $\Min$ of $\P$ or $\tP$ 
denotes the set of orbits that are minimizers.
We define  $\rho(\tilde \phi_t(z);
T_1,T_2) = \rho(\gm;T_1, T_2)$,
$\omega(z) = \gamma(\infty)$, and $\alpha(z) = \gm(-\infty)$ (when the
latter two exist).

The E-L equations are turned into a vector field using local coordinates in 
which the second coordinate is the velocity. Thus if $\pi_2$ is projection 
on the second component (\ie\ 
$z = (x,v,\tau)\in \P$, and $\pi_2(z) = v$), then for
an orbit $\phi_t(z)$, $\pi_2(\phi_t(z)) = d(p(\phi_t(z))/dt$.

The space of all the geodesics in $\Poin$ is denoted $\Geos$.
This space will always be given the topology of Hausdorff 
convergence on compact subsets. With this topology $\Geos$  is 
homeomorphic
to $S_\infty\times S_\infty - \{\hbox{diagonal}\}$ where a 
geodesic $\GM$ is identified with the pair
 $(\GM(-\infty),\GM(\infty))$. Recall that we have fixed a hyperbolic 
metric $g$ on $M$. Its geodesic flow is defined on $T_1 M$ (and is in fact 
the restriction to the invariant set $T_1 M$ of the E-L
flow of the Lagrangian $L=\norm{v}_g^2$). The geodesic flow $g_t$
 lifts to a flow $\tilde g_t$ on $T_1\Poin$.

Two flows $(X,\phi_t)$ and $(Y,\psi_t)$ are said to be 
{\de semiconjugate} (or sometimes orbit semi-equivalent) if there
is a continuous surjection $f:X\ra Y$ that takes orbits of
$\phi_t$ to those  of $\psi_t$ preserving the direction of the flow, 
but not 
necessarily the time parameterization. Note that $f$ is locally
injective when  restricted to an orbit of $\phi_t$, but $f$ may take many 
orbits of  $\phi_t$ to the same orbit of $\psi_t$. 

Given a point in $\Poin$  and a geodesic, (hyperbolic) orthogonal 
projection sends the point to a point on the geodesic. To get a image point 
in the unit tangent bundle we define $\Sigma:\Geos\times\Poin\ra 
T_1\Poin$ via
$\Sigma(\GM,z) = (x,v)$ where $x$ is  the orthogonal
projection of $z$ onto $\GM$ and $v$ is the unit vector tangent to 
$\GM$ at $x$.
\medskip
\noindent{\bf \S 3.2 Proof of Theorem B.} 
Given the Lagrangian $L$, find $K_0$ as in Proposition 2.1. Now fix
$K>K_0$ and let $\lambda, \epsilon, k'', K'',$ and $ N_0$ (all depending 
on $K$) be as in that proposition. Define the set $Q_K\subset\tP$ as the
set of $z$ that satisfy 
\medskip
(1) The orbit $\tilde \phi_t(z)$ is a minimizer and a 
$(\lambda,\epsilon)$-quasi-geodesic.

(2)  $k'' \le \rho(\tilde \phi_t(z); T_1, T_1 +  T)$ for all $T_1$,
whenever $T\ge N_0$.

(3)  $\norm{\pi_2(\tilde\phi_t(z))} \le K''$, for all $t\in\RR$.
\medskip
Note that $Q_K$ is $\tilde\phi_t$ invariant and closed and  
$\pi(Q_K)\subset\P$ is compact. 
Since each orbit in $Q_K$ is a $(\lambda, \epsilon)$-quasi-geodesic, 
by 
Theorem 1.2 there is  a constant $\kappa$  and for each $z\in Q_K$ a 
unique geodesic denoted $\Gamma_z$ with $\GM_z(\infty) = \omega(z)$,   
$\GM_z(-\infty) = \alpha(z)$, and  $d(\GM_z, p(\tilde\phi_t(z)))
\leq\kappa$.  Thus $z_i\ra z$ in $Q_K$ implies that 
$\GM_{z_i}\ra\GM_z$, 
and so the map $Q\ra\Geos$ given by
$z\mapsto\GM_z$ is continuous. 

This implies that  the ``projection'' $\sigma: Q_K\ra T_1 \tM $ defined by 
$\sigma(z) = \Sigma(\GM_z, p(z))$ is also continuous. In addition, by 
construction, $\sigma$ takes orbits of $\tilde\phi_t$ to those of the 
lift of the geodesic flow $\tilde g_t$. Also, $\sigma$ is equivariant, 
{\it 
i.e.} it descends to a map $\pi(Q_K)\ra T_1 M$. Further, by Theorem
A, $\sigma$ is onto. Unfortunately, 
$\sigma$ does not preserve the direction of time
as it is perhaps not locally injective when restricted to 
 the orbits of $\phi_t$.  This is remedied using an averaging
technique due to Fuller [F].

Fix a parameterization by arclength for each geodesics in $\Poin$. We 
will use
the parameterization to add and subtract elements on the geodesics.
Given $z\in Q_K$ and $t\in\RR$, let $a(z,t) = 
\sigma(\tilde\phi_t(z))
- \sigma(z)$, or equivalently, $a(z,t)$ is the unique $s\in\RR$ with
$\tilde g_s(\sigma(z)) = \sigma(\tp_t(z))$.  
Note that $a$ is an additive cocycle for
$\tp_t$, {\it i.e.} $a(z, t_1 + t_2) =
a(z, t_1) + a(\tilde\phi_{t_1}(z), t_2)$, for all $t_1,t_2$. 
Given ${\alpha_1} > 0$ define 
$$
\bar\sigma_{\alpha_1}(z) = \sigma(z) + 
{1\over{\alpha_1}}\int_0^{\alpha_1} a(z,t)\; dt.
$$
Equivalently, $\bar\sigma_{\alpha_1}(z) = \tilde g_s(\sigma(z))$, where 
$s={1\over{\alpha_1}}\int_0^{\alpha_1} a(z,t)\; dt$.
Informally,  $\bar\sigma_{\alpha_1}(z)$ is the average value of 
$\sigma$ over
 the orbit segment $\tilde\phi_{[0,{\alpha_1}]}(z)$. 

Now 
since for every $z\in Q_K$ we have that $\omega(z) = 
\GM_z(\infty)$, it
follows that for each $z$ there is an $\alpha_z$ so that $a(z,\alpha_z)
> 0$. Since $\pi(Q_K)$ is compact, we may find an $\alpha$ with
 $a(z, \alpha) > 0$ for all $z\in Q_K$. Let $\bs =
\bar\sigma_{\alpha}$. Now $\bs$ is clearly continuous, equivariant,
onto and takes orbits to orbits. We will show that it is injective on
orbits of $\tilde\phi_t$ by showing that for any $\beta>0$ and
$z\in Q_K$, $\bs(\tp_\beta(z))- \bs(z) > 0$.
$$
\eqalign{\bs(\tp_\beta(z))- \bs(z) &= \sigma(\tp_\beta(z))- 
\sigma(z) + 
{1\over\alpha}\left(\int_0^\alpha (a(z, \beta + t) -a(z,\beta)) - 
\int_0^\alpha a(z,t)\right)\cr
& ={1\over\alpha}\left(\int_0^\alpha a(z, \beta + t)-
\int_0^\alpha a(z,t)\right)\cr
& ={1\over\alpha}\left(\int_\beta^{\alpha+\beta} a(z,t) - 
\int_0^\alpha
a(z,t)\right)\cr
& = {1\over\alpha}\left(\int_\alpha^{\alpha+\beta} a(z,t) - 
\int_0^{\beta}
a(z,t)\right)\cr
& = {1\over\alpha}\int_0^\beta a(z, t + \alpha) - a(z,t)\cr
& = {1\over\alpha}\int_0^\beta a(\tp_t(z), \alpha)\cr
& > 0,\cr}
$$
where all integrals are with respect to $t$,
 and in the first and sixth equalities we used the cocycle equation
for $a$.

Thus for  each $K>K_0$ we have a $\pi(Q_K)$ with 
$(\pi(Q_K),\phi_t)$ 
semiconjugate to $(T_1 M, g_t)$. The set of all such $\pi(Q_K)$  
clearly 
contains a sequence $X_i$ as needed for the Theorem. \QED
\medskip
{\bf Remark:} A different perspective can be gained on Theorems A and B 
by considering the simple case $L(x, \dx, t) = \norm{\dx}^2 - V(x,t)$,
where $\norm{\cdot}$ comes from the hyperbolic metric. In this case, for
large velocities (high up in the tangent bundle),
 the time-one map  $\Phi$ of the E-L flow 
(considered as a map from $TM\ra TM$) may be thought of
as a small perturbation of the time-one map $G$ of the full
geodesic flow of the hyperbolic metric. In $TM$ above each periodic geodesic,
 $G$ has a normally hyperbolic annulus. Restricted to this annulus, $G$ is a twist map. Under small perturbation high up in the bundle,
 one expects this annulus to persist and $\Phi$ restricted to it
will also be a twist map. Using \AM\ theory on these various annuli
and taking limits one sees that
 the dynamics of $F$ near infinity to reflect those of $G$. 

 There are a number of technical problems with 
making these arguments precise, but
almost certainly these can be overcome.
 However, we prefer the techniques used here
as they are-self contained, more general and allow for fairly explicit estimates on where the persistence occurs.
\medskip

\noindent{\bf \S 3.3 Semiconjugacies, time changes and ergodic 
measures.} 
Theorems A and B give the existence of a large collection of 
minimizers, but the theorems  are somewhat unsatisfactory because 
they do not yield minimizers for which an asymptotic speed necessarily
exits. To obtain results of this type,
either one needs a great deal of control over the minimizers (for 
twist maps of the annulus this comes from the low dimensionality), or else one uses Ergodic Theory. In this subsection we briefly consider the latter.

The first ingredient, contained in Theorem B,  is the semiconjugacy 
$\bs$ from a compact $\pt$-invariant set $X_i$ onto $(T_1 M, g_t)$.
To simplify the exposition, let
us fix an $X_i$ and call it $X$, denote the semiconjugacy by $f$, and change the name of time on $X$ to $s$, so the flow restricted to $X$ is denoted $\phi_s$. A well known construction using cocycles 
(\eg\ see  [Pa] or  [HK]) allows us to perform a
time change on  the flow $\phi_s$ to obtain a new flow $\hat\phi_t$ such that $f$ is a  time-preserving semiconjugacy, \ie\ $f \hat\phi_t = g_t f$ for
all $t$. Further, to each ergodic $\phi_s$-invariant measure $\mu$ there 
corresponds an  ergodic $\hat\phi_t$-invariant measure $\hat\mu$ with
$\mu$ and $\hat\mu$ mutually absolutely continuous. 

Since $f$ is continuous, for each ergodic $g_t$-invariant $\eta$, the
set $f_*^{-1} (\eta)$ is a nonempty compact convex set in  the weak topology on measures.
The extreme points of this set are ergodic $\hat\phi_t$-invariant 
measures. Thus to each ergodic $g_t$-invariant measure $\eta$ there is at least one ergodic $\phi_t$-invariant 
$\mu$ with $f_*(\hat\mu) = \eta$. We say that $\eta$
{\de corresponds} to $\mu$.

The second ingredient  is a way of measuring the progress of orbits
in the universal cover.
Given  $z\in \P$ and $t\in \reals$, pick a lift $\tz\in\tM$ and let
$D(z,t) = d(p(\tpt(\tz)), p(\tz))$. Note that the definition of $D$ is independent of the choice of lift $\tz$ and that $D(z, t) = \rho(z; 0, t)$, 
with $\rho$ as defined above. Now let 
$$D^\ast(z) = \lim_{t\ra\infty} {D(z,t)\over t}$$
if the limit exists.
The triangle inequality for the metric $d$ implies that
$D$ satisfies $D(z, t + s) \leq D(z, t) + D(\phi_t(x), s)$
for all $z$, $s$ and $t$. Thus $D$ is a subadditive cocycle for 
$\phi_t$, and so by Kingman's subadditive ergodic theorem 
(see, for example, [Po]) we have that $D^\ast$ exists and is constant almost everywhere 
with respect to an ergodic $\phi_t$-invariant measure $\mu$.

Recalling the original situation from Theorem B, we now 
see that for each $i$ and for each  ergodic $g_t$-invariant measure
 $\eta$ there
corresponds a $\phi_t$-invariant measure $\mu$ supported on $X_i$ such that $D^\ast(z)$ exists and is constant for $\mu$ almost every point $z$. Further, using the speed bounds from the theorem, we have $k_i \leq D^\ast(z) \leq
k_{i+1}$.

Given the ergodic measure $\mu$,  the pair
$(D^\ast(\mu), \eta)$, where $\eta$  corresponds to $\mu$, can be interpreted 
as giving the length and direction of a kind of rotation vector. 
One can see from the definition that the  correspondence between ergodic
measures invariant under $\phi_t$ and $g_t$ means vaguely that the dynamical behavior of one ``shadows'' the other.
What is lacking is a precise
meaning of this correspondence which doesn't require
 {\it a priori} knowledge of a semiconjugacy.
This lack is remedied  in a subsequent paper ([Bd]).

\bigskip
\centerline{\bf References}
\medskip

\bibitem{AM}
Abraham, R. and Marsden, J.E. (1985) {\it Foundations of Mechanics},
Addison-Wesley

\bibitem{AD}
Aubry, S. and Le Daeron, P.Y. (1983)  The discrete Frenkel-Kontorova
model and its extensions I: Exact results for the ground states,
 {\it Physica} {\bf 8D}, pp. 381-422

\bibitem{B1}
Bangert, V. (1988) Mather  sets for twist maps and geodesics on tori,
{\it Dynamics Reported} {\bf 1} (eds. Kirchgraber, U. and Walther, H.O.) John Wiley

\bibitem{B2}
Bangert, V.  (1989)  Minimal geodesics, {\it Ergod. Th. \& Dynam. 
Sys.},  {\bf 10}, pp. 263-286

\bibitem{BKS}  Bedford, T., Keane, M., and  Series, C. (eds) (1991)
  {\it Ergodic theory, symbolic dynamics, and hyperbolic
spaces}, Oxford University Press

\bibitem{Bn}  Benci, V. (1986) Periodic solutions for Lagrangian systems on
compact manifolds, {\it J. Diff. Eq.}, {\bf  63}, 135-161

\bibitem{BK}
Bernstein, D. and  Katok, A.  (1987) Birkhoff periodic orbits for small
perturbations of completely integrable Hamiltonian systems, 
{\it Invent. Math.} {\bf 88},  pp. 225-241

\bibitem{Bd}
Boyland, P, Asymptotic dynamical invariants on hyperbolic manifolds,
 in preparation

\bibitem{BG} Boyland, P. and Gol\'e, C. (1995) Dynamical stability in 
Lagrangian systems, Proceedings of the NATO Advanced Study Institute on 
Hamiltonian systems with three or more degrees or freedom,  S'Agar\'o,
Spain (to appear)
 
\bibitem{BP}
Bialy, M. and Polterovitch, L. (1992) Hamiltonian systems, Lagrangian tori,
and\hfill\break Birkhoff's theorem {\it Math. Ann.} {\bf 292}, pp. 619-627

\item{[CDP]} Coornaert, M., Delzant, T. and Papadopoulos, A. (1990)    
    {\it G\'eom\'etrie et theorie de groupes: les groups hyperboliques de Gromov}, Springer-Verlag

\bibitem{DM}
Denvir, J. and MacKay, R.S. (1995) Consequences of contractible
geodesics on surfaces, preprint

\bibitem{F} 
Fuller, F. B. (1965) On the surface of section and periodic trajectories,
{\it Am. J. Math.}, {\bf 87}, pp. 473--480

\bibitem{GH}
Ghys, E. and de la Harpe, P. (eds) (1990)  {\it Sur les groupes hyperboliques d'apr\`es Mikhael Gromov}, Birkh\"auser

\bibitem{Gl}
Gol\'e, C. (1994) Periodic orbits for Hamiltonians in cotangent
bundles {\it Trans. A.M.S.} {\bf 343},  Number 1, pp.  327-347 

\bibitem{G1} 
Gromov, M. (1977) Three remarks on geodesic dynamics
and fundamental group,\hfill\break preprint, SUNY at Stony Brook

\bibitem{G2} Gromov, M., Hyperbolic groups (1987)  {\it Essays in Group Theory}, ed. by Gersten, S. M., MSRI Publications no. 8, Springer-Verlag

\bibitem{HK} Hasselblat, B. and Katok, A. (1995) {\it Introduction to
the
modern theory of dynamical systems}, Cambridge University Press.

\bibitem{K}
Katok, A.  (1992) Minimal orbits for small perturbations
of completely integrable Hamiltonian systems, 
{\it Twist mappings and their applications, I.M.A.} {\bf 44} (McGehee, R.
and Meyer, K. eds) Springer-Verlag

\bibitem{MMS}
MacKay, R.S., Meiss, J.D. and Stark, J. (1989) Converse KAM theory for symplectic twist maps, {\it Nonlinearity} {\bf 2}, pp. 555-570

\item{[Mn]}
Ma\~n\'e, R. (1991) {\it Global variational methods in conservative dynamics},  $18^o$ Col\'oquio Brasileiro de Matem\'atica, IMPA

\item{[Ma1]}
Mather, J. (1991)  Action minimizing invariant measures for
positive definite Lagrangian systems, {\it Math. Z.} {\bf 207}, pp. 169-207

\item{[Ma2]}
Mather, J. (1993)  Variational constructions of connecting
orbits, {\it Ann. Inst. Fourier, Grenoble} {\bf 43}, 5, pp. 1349-1386

\bibitem{Me}
Meiss, J. D. (1992) Symplectic maps, variational principles, and transport, {\it Rev. Mod. Phys.} {\bf 64}, No. 3, pp. 795-848

\bibitem{M}
Morse, M. (1924) A fundamental class of geodesics on any closed
surface of genus greater than one, {\it Trans. Math. Soc.} {\bf 26}, pp. 25-60

\bibitem{Mo} Moser, J.  (1986) Monotone twist mappings and the calculus 
of variations, {\it Erg. Th. \& Dynam. Sys. } {\bf 8*}, pp. 199-214

\bibitem{Pa} Parry, W. (1981)  {\it Topics in ergodic theory}, Cambridge
University Press.

\bibitem{Po} Pollicott, M. (1993)
{\it Lectures on ergodic theory and Pesin theory on compact manifolds},
Lond. Math. Soc. Lecture Notes in Math, {\bf 180}, Cambridge University
Press.

\bigskip
\medskip

\noindent {\bf Phililp Boyland},
Dept. of Theoretical
and Applied Mechanics, 
University of Illinois, 
Urbana, IL 61801 \ \ 
boyland@math.uiuc.edu
\medskip
\noindent {\bf Christophe Gol\'e},
Dept. of Mathematics,
University of California,
Santa Cruz, CA 95064 \ \ 
gole@cats.ucsc.edu

\bye